\documentclass[review, 12pt]{elsarticle}

\usepackage{algpseudocode}
\usepackage[plain]{algorithm}
\usepackage{amssymb, amsmath, amsthm, amscd, array}
\usepackage{balance, booktabs, bm}
\usepackage[labelsep=period]{caption}
\usepackage{calc, color, comment}
\usepackage{cases}
\usepackage{empheq, enumitem, epsfig, epstopdf, etoolbox}
\usepackage{flafter}
\usepackage{float}
\usepackage{geometry}
\usepackage{graphicx}
\usepackage[colorlinks=true, linkcolor=blue, citecolor=blue]{hyperref}
\usepackage[utf8]{inputenc}
\usepackage{latexsym, lipsum, lineno}
\usepackage{mathtools}
\usepackage{makecell, multirow}
\usepackage{natbib}
\usepackage{pdflscape}
\usepackage{pict2e}
\usepackage{rotating}
\usepackage{setspace, subcaption}
\usepackage{tabularx}
\usepackage[flushleft]{threeparttable}
\usepackage{wrapfig}
\usepackage[dvipsnames]{xcolor}
\usepackage{xtab}
\usepackage{colortbl}
\makeatletter
\renewcommand{\ALG@name}{Scheme}
\makeatother

\makeatletter
\renewcommand{\ALG@beginalgorithmic}{\footnotesize}
\makeatother

\makeatletter
\newcommand{\multiline}[1]{%
  \begin{tabularx}{\dimexpr\linewidth-\ALG@thistlm}[t]{@{}X@{}}
    #1
  \end{tabularx}
}
\makeatother

\numberwithin{equation}{section}

\modulolinenumbers[1]

\begin{document}

\begin{frontmatter}

\title{A Stochastic Differential Equation Model for Predator-Avoidance Fish Schooling}

\author[aff1]{Aditya Dewanto Hartono}

\author[aff3]{Linh Thi Hoai Nguyen}

\author[aff1,aff2]{T\^{o}n Vi\d {\^{e}}t T\d{a}}

\address[aff1]{
Mathematical Modeling Laboratory, Department of Agro-environmental Sciences, Kyushu University
}

\address[aff2]{
Center for Promotion of International Education and Research, Kyushu University
}

\address[aff3]{
Institute of Mathematics for Industry, Kyushu University
}

\begin{abstract}
This paper presents a system of stochastic differential equations (SDEs) as mathematical model to describe the spatial-temporal dynamics of predator-prey system in an artificial aquatic environment with schooling behavior imposed upon the associated prey. The proposed model follows the particle-like approach where interactions among the associated units are manifested through combination of attractive and repulsive forces analogous to the ones occurred in molecular physics. Two hunting tactics of the predator are proposed and integrated into the general model, namely the center-attacking and the nearest-attacking strategy. Emphasis is placed upon demonstrating the capacity of the proposed model in: (i) discovering the predator-avoidance patterns of the schooling prey, and (ii) showing the benefit of constituting large prey school in better escaping the predator's attack. Based on numerical simulations upon the proposed model, four predator-avoidance patterns of the schooling prey are discovered, namely Split and Reunion, Split and Separate into Two Groups, Scattered, and Maintain Formation and Distance. The proposed model also successfully demonstrates the benefit of constituting large group of schooling prey in mitigating predation risk. Such findings are in agreement with real-life observations of the natural aquatic ecosystem, hence confirming the validity and exactitude of the proposed model.
\vspace*{0.2cm}
\end{abstract}

\begin{keyword}
Predator-Prey System \sep Swarm Behavior \sep Particle-Based Model \sep Stochastic Differential Equations \sep Predator-Avoiding Patterns \sep Fish Schooling 
\vspace*{0.2cm}
\MSC[2020] 92-10 \sep 60H10 \sep 68W10
\end{keyword}

\end{frontmatter}

\allowdisplaybreaks

\vspace*{0.1cm}
\section{Introduction} \label{introduction}
\noindent Swarm dynamics, one of the most commonly observed phenomena in the real world, have attracted the interest of many researchers from diverse fields including biology, mathematics, and computer engineering \citep{Aoki1982School, Reynolds1987Flocks, Cucker2007Mathematics, Couzin2009Collective, Vicsek2012CollectiveMotion, TonandLinh2018Foraging}. As one manifestation of collective behavior in nature, swarm dynamics offer many advantages to the associated individuals: higher foraging success \citep{Pitcher1982FishFood, Jolles2020Role, Ioannou2021Grouping}, improve protection upon predators \citep{Jolles2020Role, Rubenstein1978Group, Partridge1982Structure, Pitcher1986Functions, Magurran1990Adaptive, Demsar2015PredatorAttacks}, increase mating chance \citep{Jolles2020Role, Rubenstein1978Group}, and more efficient energy consumption due to hydrodynamic or aerodynamic advantages \citep{Jolles2020Role, Ioannou2021Grouping, Rubenstein1978Group}. Based on observational and empirical investigations of the interaction between congregating animals and their neighbors, researchers have constructed artificial life \citep{Kim2006Comprehensive}, for example, mathematical models of schooling prey in the presence of a predator.  

As far as mathematical modeling of predator-prey interaction with schooling behavior upon the associated prey is concerned, two different approaches of modeling are routinely considered. The first approach examines the predator-prey system from the standpoint of tracking population densities of the associated entities. In this framework, mathematical models are typically constructed by integrating functional responses into the general description (see, for example, \citep{Abrams2000Evolution, Abrams2000Nature, Bera2015Modelling, Bera2016Dynamics, Maiti2016Predator, Maiti2016Herd, Manna2018Analysis, Manna2019Deterministic}). The second approach of modeling examines the predator-prey system as self-propelled units in an artificial ecosystem. Here, the associated entities are treated as agents (or particles) that move around and interact with each other perpetually in the prescribed domain. Hereinafter, the discussion concerning mathematical modeling of the predator-prey system follows the latter approach.

Let us review some existing studies. \citet{Oboshi2002Collective} proposed a mathematical model to describe predator-prey interaction in an aquatic ecosystem. Therein, they designed an artificial habitat where schooling prey fish coexist with solitary predator. They used a local rule that each prey fish in the school chooses one way of action among four possibilities according to the closest mate to construct a model of difference equations for anti-predator schooling behavior. Their model is able to return single predator-avoidance maneuver that fits real observation at the Port of Nagoya Public Aquarium. Such an evasive maneuver is also observed in \citep{Pitcher1986Functions, Vabo1997Individual, Capuska2011Dive, Parrish2002FishSchools, Hulsman2021FishyTales}. Later on, the same authors \citep{Oboshi2003Simulation} extended their assessment by incorporating three different conditions for prey-predator encounters in the artificial habitat, namely predator faces center, side, and vicinity of the schooling prey fish. By assimilating such configurations into their model, they managed to obtain two additional anti-predator patterns.

\citet{Nishimura2000Studying} investigated mechanisms for a predator to select an individual prey from a group. Therein, the author introduced a priority function as pertinent parameter determining the selectivity of predator upon a particular prey in the group. The study is later extended \citep{Nishimura2002PredatorSelection} by including three different priority functions, each corresponds to the selection of nearest, peripheral, and split individuals among the schooling prey. It is shown that, from the predator's viewpoint, the selection of peripheral individual is the best of the three available options: such a strategy yields the highest average number of successfully captured prey.

\citet{Zheng2005Behavior} considered the evasive maneuver of schooling prey fish as complex collective traits resulting from the combination of schooling, cooperative, and selfish behavior. Three mathematical models of prey fish behavior are constructed, each corresponds to a fundamental trait. A combination of such models yields one anti-predator pattern which is similar to the one obtained by \citet{Oboshi2002Collective}.

\citet{Lee2006Dynamics} implemented molecular dynamics (MD) to model the dynamical behavior of schooling prey in response to predator's attack. They retrieved one predator-avoidance pattern of the schooling prey, namely the crescent-shaped maneuver. In a separate study \citep{Lee2006Predator}, the transitional regimes associated with such a maneuver are elaborated further.

\citet{Demvsar2014Simulated} studied three hunting tactics of the predator (attack center, nearest, and isolated prey) in association with two distinct behavior of the prey (schooling and solitary) using a mathematical model based on fuzzy logic. They found that compared to individualistic behavior, schooling is the optimal defense mechanism to avoid predation. Later on, \citet{Demsar2015PredatorAttacks} extended the assessment by including composite hunting tactics of the predator. Two composite hunting tactics were proposed: (i) the predator may choose one of three simple hunting strategies (attack nearest, center, and peripheral prey) in successive attacks based on probability, and (ii) the predator was set to initially disperse the schooling prey and then pursue an isolated individual. They concluded that predator's confusion plays an important role upon the evolution of the composite tactics. Moreover, they conveyed that when confusion is taken into account, the latter composite tactic is the better predation mode that provides favorable outcome for the predator.

Despite the vibrant research activity in the field of mathematical modeling of predator-prey system, up to our knowledge, there is still no mathematical model that provides not only one but all real predator-avoidance patterns mentioned above. Moreover, there is also no mathematical model which shows the benefit of making a large school of fish in surviving predator's attack. In this paper, we aim to construct a mathematical model based on stochastic differential equations (SDEs) for predator-avoidance fish schooling which can show not only simulated patterns adapting to the real patterns but also the benefit of making a large school of fish in surviving predator's attack. The SDE model is chosen as it is well-suited for capturing complex movements and interaction rules of individuals in diverse environment. Furthermore, an SDE model is amenable for undertaking numerical simulations in computer platform.

In \citep{Uchitane2012ODE, Linh2015Quantitative}, we introduced a system of stochastic differential equations to describe the process of fish schooling in a free, unbounded environment based on a set of behavioral rules of fish schooling proposed by \citet{Camazine2001Biological}. The rules are stated as follows:
\vspace*{-0.2cm}
\begin{enumerate}[label=(\alph*)]
\itemsep0cm 
\item The school has no leaders and each fish follows the same behavioral rules.
\item Each fish uses some form of a weighted average of the position and orientation of its nearest neighbors to decide where to move.
\item There is a degree of uncertainty in the individual's behavior that reflects both the imperfect information-gathering ability of a fish, and the imperfect execution of the fish's actions.
\end{enumerate}
\vspace*{-0.2cm}
However, these rules do not consider obstacle-avoiding mechanisms and the foraging process of a school of fish. In accordance with this, in \citep{TonandLinh2018Foraging, Linh2016ObstacleAvoiding}, we proposed behavioral rules for both phenomena, respectively. Furthermore, using the model in \citep{Uchitane2012ODE}, we constructed two mathematical models for obstacle-avoiding mechanisms and the foraging behavior of fish schools. As a result, four obstacle-avoiding patterns are discovered. In addition, the probability of foraging success is observed to escalate up to an optimal value of school size, then gradually decreases as the school size gets larger.

In order to construct the mathematical model for predator-avoidance fish schooling in the present paper, we integrate two hunting strategies of a predator into our model in \citep{Uchitane2012ODE}. The first strategy is that the predator always moves to the center of the prey fish school and catches the neighboring fish while moving. The second one is that it chases the \textquotedblleft nearest\textquotedblright\ prey fish. For the latter strategy, the predator chases all prey fish in the group with imposed weighted coefficients that are inversely proportional to the distance from it to the corresponding individuals. Therefore, in this study, the \textquotedblleft nearest\textquotedblright\ attacking strategy means that the predator pursues all the prey fish in the school while giving utmost attention to its nearest prey.

The paper is organized as follows. In Section \ref{model_description}, we introduce our SDE model for predator-avoidance fish schooling and explain how we construct mathematical expressions for the two hunting strategies and integrate them into the general model. We also explain \textquotedblleft being eaten\textquotedblright\ conditions for the prey fish. Section \ref{predator_avoiding_patterns} provides four predator-avoidance patterns obtained from the model which fit real-life anti-predatory maneuvers. We call them Split and Reunion, Split and Separate into Two Groups, Scattered, and Maintain Formation and Distance. In addition, pseudocodes to obtain these patterns are also given. Section \ref{benefit_schooling} demonstrates the capability of our model in describing the benefit of making a large school of fish for better surviving the predator's attack\textemdash an advantage that is particularly well-recognized in real-life situation\textemdash by means of probability of survival of prey fish and the average time of prey being alive when the school size increases. Finally, Section \ref{conclusions} poses the concluding remarks of the present study.

\section{Model Description} \label{model_description}
\noindent In this section, we newly introduce an SDE model to describe the dynamics of schooling prey under the attack of a predator in an obstacle-free domain $\mathbb{R}^{d} \ (d = 2 \ \textrm{or} \ 3)$. In constructing the model, we adopt the particle-like approach; namely, each of the associated units in the predator-prey system is treated as if they were particles moving in space $\mathbb{R}^{d}$ where interactions among them are manifested in the forms of fundamental forces analogous to the ones occurred in molecular physics \citep{Vicsek2012CollectiveMotion, Linh2021ReviewModels}. The proposed model is formulated as follows:\\
\begin{equation} \label{eq2-1}
\begin{aligned}
\begin{cases}
d{x}_{i}(t) =& v_{i}dt + \sigma_{i}d{w}_{i}(t), \quad i = 1, 2, \dots, N, \\
d{v}_{i}(t) =&\Bigg[
-\alpha \sum \limits_{j = 1, j \ne i}^{N} \left(\frac{r^{p}}{\|x_{i} - x_{j}\|^{p}} - \frac{r^{q}}{\|x_{i} - x_{j}\|^{q}}\right) \left(x_{i} - x_{j}\right) \\
&-\beta \sum \limits_{j = 1, j \ne i}^{N} \left(\frac{r^{p}}{\|x_{i} - x_{j}\|^{p}} + \frac{r^{q}}{\|x_{i} - x_{j}\|^{q}}\right) \left(v_{i} - v_{j}\right)  \\
&+H\left(x_{i}, y\right)
\Bigg]dt, \quad i = 1, 2, \dots, N, \\
dy(t) =& vdt + \sigma d{w}_{t}, \\
dv(t) =& F\left(x_{i}, v_{i}, y, v\right)dt. 
\end{cases}
\end{aligned}
\end{equation}

\noindent Here, $N$ depicts the population of the schooling prey; $x_{i}(t)$ and $v_{i}(t)$ $(i = 1, 2, \dots, N)$ respectively denote the position and velocity of the $i$-th prey at time $t$; $y(t)$ and $v(t)$ correspondingly represent the position and velocity of the predator at time $t$, while $\|\cdot\|$ specifies the Euclidean norm of a vector. 

As can be seen in the system \eqref{eq2-1}, the proposed model consists of two parts: (i) a group of expressions that describe the dynamics of the schooling prey (the first two equations) and (ii) a group of expressions that depict the behavior of the predator (the last two equations). Let us discuss each of these groups of the model in more detail.

\subsection{The Schooling Prey} \label{model_prey}
\noindent The dynamics of the schooling prey are described by the first two equations in \eqref{eq2-1}. These expressions are extended versions of the ones developed in our previous study \citep{Uchitane2012ODE}. The first equation is an SDE that tracks the position $x_{i}(t)$ of the $i$-th prey fish. It contains a stochastic term $\sigma_{i}d{w}_{i}(t)$, considered as a white noise that reflects uncertainty in an individual's behavior within the school due to some kind of defect in the process of information-gathering among the individuals and imperfect execution of an individual's action. 

Subsequently, the second equation in  \eqref{eq2-1} solves the velocity $v_{i}(t)$ of the $i$-th prey fish. Therein, $\alpha$ and $\beta$ are positive coefficients; parameters $1 < p < q < \infty$ denote positive exponents; and $r > 0$ is the critical distance between two individuals in the school. Together with the Euclidean distance $\|x_{i} - x_{j}\|$, the overall term $-\Big(\frac{r^{p}}{\|x_{i} - x_{j}\|^{p}} - \frac{r^{q}}{\|x_{i} - x_{j}\|^{q}}\Big) \left(x_{i} - x_{j} \right)$ determines interaction between two prey individuals $i$ and $j$ in the form of attraction or repulsion. If the distance exceeds $r$, namely $\|x_{i} - x_{j}\| > r$, then they would move towards each other due to the induced attractive force between them. Contrarily, if $\|x_{i} - x_{j}\| \leqslant r$, then the interaction would be the repulsive force so as to avoid collision.

To account for the effect of the predator's presence upon the behavior of the schooling prey, the term $H (x_{i}, y)$ is introduced within the expression of $v_{i}(t)$. Overall, such a remark manifests repulsive interaction between each individual prey and the predator. In this study, we propose the following mathematical formula for $H(x_{i}, y)$:\\
\begin{equation} \label{eq2-5}
H(x_{i}, y) = \delta \frac{R_{1}^{\theta_{1}}}{\|x_{i} - y\|^{\theta_{1}}} (x_{i} - y).
\end{equation}

\noindent Here, $R_{1} > r$ is a constant that reflects the threshold distance where an individual prey begin to display the anti-predation maneuver to avoid the approaching predator; $\delta$ and $\theta_{1}$ are positive coefficients. The strength of repulsion between the schooling prey and the predator depends upon the distance between them, $\|x_{i} - y\|$. The closer the distance, the stronger the induced repulsive force between them. 

\subsection{The Predator} \label{model_predator}
\noindent Let us now turn to discuss the expressions that represent the behavior of the predator in our model. Analogous to the prey counterpart, the position $y(t)$ of the predator is described by the third equation in \eqref{eq2-1}. Therein, the term $\sigma d{w}_{t}$ represents uncertainty in predator's action due to the imperfect decision-making during hunting activity.

The velocity of predator $v(t)$ in \eqref{eq2-1} is formulated based on generalized Newton's second law of motion. Therein, the resultant force term $F(x_{i}, v_{i}, y, v)$ can be regarded as an attractive driving force that draws the predator towards the schooling prey. Apart from its conventional definition in Newton's second law of motion, in our model, the force term $F(x_{i}, v_{i}, y, v)$ is specifically designed to mimic the hunting strategy adopted by the predator. In this regard, we propose two hunting tactics for the predator, namely (i) the center-attacking and (ii) the nearest-attacking strategy. Each of these tactics occupied its own mathematical formula for $F(x_{i}, v_{i}, y, v)$.

\subsubsection{Hunting Strategy I: Predator Attacks the Center of the Schooling Prey} \label{predator_center_attacking}
\noindent In an aquatic ecosystem, this hunting strategy is mainly exemplified by swordfish (\textit{Xiphias gladius}) \citep{Demsar2015PredatorAttacks, Pavlov2000PatternsSchooling}. Observation \citep{Pavlov2000PatternsSchooling} revealed that a similar hunting tactic is also adopted by other predatory fish species such as tuna (\textit{Thunnus} spp.), barracuda (\textit{Sphyraena} spp.), marlin (\textit{Makaira} spp.), and sailfish (\textit{Istiophorus platypterus}). 

The associated mathematical expression of $F(x_{i}, v_{i}, y, v)$ for this hunting strategy is devised as:\\
\begin{equation} \label{eq2-6}
F(x_{i}, v_{i}, y, v) = -\frac{R_{2}^{\theta_{2}}}{\|y - x_{c}\|^{\theta_{2}}} 
\left[\gamma_{1}(y - x_{c}) + \gamma_{1} \gamma_{2} (v - v_{c})\right]. 
\end{equation}

\noindent Here, $R_{2} > r > 0$ reflects the threshold distance where the predator begins to display predation maneuver to capture the nearby prey; $\theta_{2}$, $\gamma_{1}$, and $\gamma_{2}$ are positive coefficients; the center position $x_{c}$ and velocity $v_{c}$ of the schooling prey are defined respectively as the average positions and velocities of the individuals constituting the school:\\
\begin{equation} \label{eq2-7}
x_{c} = \frac{1}{N} \sum \limits_{i = 1}^{N} x_{i}, \quad
v_{c} = \frac{1}{N} \sum \limits_{i = 1}^{N} v_{i}.
\end{equation}

\noindent The strength of attractive force depends upon the distance between the predator and the central positional point of the schooling prey, $\|y - x_{c}\|$. The closer the distance between them, the stronger the triggered attractive force that draws the predator towards the center of the schooling prey. 

\subsubsection{Hunting Strategy II: Predator Attacks the Nearest Prey} \label{predator_nearest_attacking}
\noindent The second hunting strategy is devised in such a way that allows the predator to pay more attention to the nearest (visible) prey. Put in other words, the predator focuses to hunt the most vulnerable target from the schooling prey. In the marine ecosystem, this particular predation tactic is mainly exhibited by the camouflage-fish species such as stonefish (\textit{Synanceia verrucosa}) \citep{Grobecker1983Feed}, scorpionfish (\textit{Iracundus signifier}) \citep{Shallenberger1973Luring}, frogfish (\textit{Antennarius commerson}) \citep{Pietsch2020Frogfishes}, and the Atlantic stargazer fish (\textit{Uranoscopus scaber}) \citep{Rizkalla2008Feeding}.

For this particular hunting strategy, the corresponding mathematical expression for $F(x_{i}, v_{i}, y, v)$ is given as:\\
\begin{equation} \label{eq2-8}
F(x_{i}, v_{i}, y, v) = - \frac{1}{N} \sum \limits_{j = 1}^{N} \frac{R_{2}^{\theta_{2}}}{\|y - x_{j}\|^{\theta_{2}}}
\left[\gamma_{1}(y - x_{j}) + \gamma_{1}\gamma_{2}(v - v_{j})\right]. 
\end{equation}

\noindent In (\ref{eq2-8}), all the pertinent parameters occupy similar definitions as the ones given for (\ref{eq2-6}). The striking discrepancy with the former hunting strategy, however, lies in the fact that in the current hunting tactic the predator's behavior is affected by the movements of all individual prey constituting the school. The tendency of the predator to put special attention upon the nearest prey is depicted by the expression $\frac{1}{\|y - x_{j}\|^{\theta_{2}}}$. Such a remark serves as the weighted coefficient for each of the associated prey. It implies that the closest prey would trigger the strongest attraction to the predator.

\subsection{\textquotedblleft Being Eaten\textquotedblright\ Condition for the Prey} \label{eaten_condition}
\noindent Having discussed the proposed model, this section proceeds to elaborate the condition when one or more prey is captured by the predator during its attack. We use the term \textquotedblleft being eaten\textquotedblright\ to describe such a situation. Mathematically, the condition of the $i$-th prey \textquotedblleft being eaten\textquotedblright\ is defined as:\\
\begin{equation} \label{eq2-9}
\|y - x_{i}\| < mr,
\end{equation}

\noindent where $0.01 \leqslant m \leqslant 1$ specifies the multiplier constant. Put it in words, a particular prey is declared to be eaten by the predator when its Euclidean distance with the approaching predator is less than the critical distance between each individual prey in a school $r$ multiplied by a positive constant $m$. From this definition, it is therefore possible for the predator to capture multiple prey simultaneously. Such a situation can be specified in the model by tuning the multiplier constant $m$.

\section{Predator-Avoiding Patterns} \label{predator_avoiding_patterns}
\noindent As outlined in the Introduction, one of our aims in this study is to implement the proposed predator-prey model to discover the predator-avoiding patterns exhibited by the schooling prey. To fulfill such an objective, we carry out numerical simulations upon our predator-prey model. The associated computer codes are constructed in the form of programming functions thereby allowing us to better organize the assessment.

In the beginning, numerical simulations are carried out to establish schooling condition for the associated prey. Each of the prey is initially set to occupy a random position in the domain $\mathbb{R}^{d}$. The first task is therefore to bring all these scattered prey together in a manner that follows certain behavioral rules such that they manage to form a school. In the light of such an objective, the predator is (deliberately) excluded from the system. For the purpose of generating schooling prey, we adopt the model from our previous studies \citep{Uchitane2012ODE, Linh2015Quantitative}:\\
\begin{equation} \label{eq3-1}
\begin{aligned}
\begin{cases}
d{x}_{i}(t) = &v_{i}dt + \sigma_{i}d{w}_{i}(t), \quad i = 1, 2, \dots, N, \\
d{v}_{i}(t) =& \Big[
-\alpha \sum \limits_{j = 1, j \ne i}^{N} \left(\frac{r^{p}}{\|x_{i} - x_{j}\|^{p}} - \frac{r^{q}}{\|x_{i} - x_{j}\|^{q}} \right)(x_{i} - x_{j})\\
&-\beta \sum \limits_{j = 1, j \ne i}^{N} \left(\frac{r^{p}}{\|x_{i} - x_{j}\|^{p}} + \frac{r^{q}}{\|x_{i} - x_{j}\|^{q}} \right)(v_{i} - v_{j})\\
&-k v_{i}
\Big]dt, \quad i = 1, 2, \dots, N.
\end{cases}
\end{aligned}
\end{equation}

This system is similar in form to the first two equations in \eqref{eq2-1} except for the last term in the velocity expression. Therein, the friction force $-kv_{i}$ supersedes the term depicting contribution from the predator. Parameter $k$ denotes a constant that manifests the magnitude of such a force. All the other quantities in \eqref{eq3-1} occupy similar definitions as the ones given for the first two equations in \eqref{eq2-1}. Scheme~\ref{algorithm_1} outlines the pseudocode for undertaking simulation to generate the schooling prey. Therein, $\sigma_{w}$ denotes the magnitude of the white noise; $\Delta t$ is the time step; $t_{\textrm{max}}$ is the maximum simulation time; and $v_{\textrm{max}}$ is the maximum allowable velocity of individual prey within the school.

\begin{algorithm}[!b]
\vspace*{0.2cm}
\caption{Pseudocode to generate the schooling prey.}\label{algorithm_1}
\setstretch{1.1}
\hrulefill
\begin{algorithmic}[1]
\Function{Schooling Generation}{$N, d, \sigma_{w}, \alpha, \beta, p, q, r, \Delta t, t_{\textrm{max}}, v_{\textrm{max}}$}
\State Initialize $x_{i}$ randomly for $i = 1, 2, \dots, N$; \Comment{\parbox[t]{0.4\textwidth}{\textit{Specify initial position for every prey}}}
\State \multiline{Initialize $v_{i}$ either randomly or static \Comment{\parbox[t]{0.4\textwidth}{\textit{Specify initial velocity for every prey}}} \\ for $i = 1, 2, \dots, N$;}
\For{$t = 1$ \textbf{to} $t = t_{\textrm{max}}$} \Comment{\parbox[t]{0.62\textwidth}{\textit{Execute the calculations in parallel manner for $i = 1, 2, \dots, N$}}}
\State Calculate random noise $\sigma_{i} = \sigma_{w} \times \sqrt{\Delta t} \times \textrm{random number of size} \ (N \times d)$;
\State Calculate $x_{i}$ using the first equation in \eqref{eq3-1};
\State Calculate $v_{i}$ using the second equation in \eqref{eq3-1};
\If{$\|v_{i}\| > v_{\textrm{max}}$} \Comment{\parbox[t]{0.4\textwidth}{\textit{Condition to avoid explosion}}}
\State Update $v_{i} \coloneqq v_{i} \times \dfrac{v_{\textrm{max}}}{\|v_{i}\|}$;
\EndIf
\EndFor
\State \multiline{Transform $x_{i}$ such that the center point lies at coordinate $(0, 0)$ (for $d = 2$)\\or $(0, 0, 0)$ (for $d = 3$);}
\textbf{return} $x_{i}$, $v_{i}$
\EndFunction
\end{algorithmic}
\hrulefill
\end{algorithm}

Subsequently, the outcomes of the schooling prey simulation in terms of prey's position $x_{i}$ and velocity $v_{i}$ are passed on as initial condition to the predator-prey function. By doing so, the predator-prey simulation can be carried out under the condition of the associated prey already formed a school. Scheme~\ref{algorithm_2} summarizes the sequential steps taken for executing the predator-prey simulation to acquire the anti-predator maneuvers of the schooling prey.

\begin{algorithm}[!b]
\vspace*{0.2cm}
\caption{\centering Pseudocode for executing predator-prey simulation\newline to capture predator-avoiding patterns.}\label{algorithm_2}
\setstretch{1.1}
\hrulefill
\begin{algorithmic}[1]
\Function{Pattern Generation}{$N, d, \sigma_{w}, \alpha, \beta, \delta, r, p, q, \theta_{1}, \theta_{2}, \gamma_{1}, \gamma_{2}, R_{1}, R_{2}, m, \Delta t, t_{\textrm{max}}$}
\State Load $x_{i}$ and $v_{i}$ $(i = 1, 2, \dots, N)$ from the output of \Call{Schooling Generation}{} function;
\State Initialize $y$ randomly;
\State Initialize $v$ randomly or static;
\State Define $N_{\textrm{survived}}(t) \coloneqq \textrm{number of survived prey at time step}\ t$;
\For{$t = 1$ \textbf{to} $t = t_{\textrm{max}}$} \Comment{\parbox[t]{0.62\textwidth}{\textit{Execute the calculations in parallel manner for $i = 1, 2, \dots, N$}}}
\State \multiline{Calculate random noise for the schooling prey:\\ $\sigma_{i} \coloneqq (\sigma_{w})_{\textrm{prey}} \times \sqrt{\Delta t} \times \textrm{random number of size} \ (N \times d)$;}
\State \multiline{Calculate random noise for the predator:\\ $\sigma \coloneqq (\sigma_{w})_{\textrm{predator}} \times \sqrt{\Delta t} \times \textrm{random number of size} \ (1 \times d)$;} 
\State Calculate $x_{i}$ using  the first equation in \eqref{eq2-1};
\State Calculate $y$ using   the third equation in  \eqref{eq2-1};
\State Calculate $v_{i}$ using  the second equation in  \eqref{eq2-1}; 
\State \multiline{Calculate $F_{i}$ using Eq.~(\ref{eq2-6}) (for hunting tactic I) or Eq.~(\ref{eq2-8}) (for hunting tactic II);}   
\State Calculate $v$ using the last equation in  \eqref{eq2-1};  
\State Apply \textquotedblleft being eaten\textquotedblright\ condition (Eq.~(\ref{eq2-9})); 
\State Calculate $N_{\textrm{survived}}(t)$ based on the above \textquotedblleft being eaten\textquotedblright\ condition;
\If{$N_{\textrm{survived}}(t) = 0$}
\State \textbf{break}; 
\Else
\State Update $x_{i}$ by including only the survived prey;
\State Update $v_{i}$ by including only the survived prey;
\EndIf
\EndFor\\
\textbf{return} $x_{i}$, $v_{i}$, $y$, $v$, $N_{\textrm{survived}}(t)$
\EndFunction
\end{algorithmic}
\hrulefill
\end{algorithm}

In the computer program, we do not implement the conventional loop technique to carry out the successive calculations. Instead, we use the vectorized adapted parallel computation to accelerate the calculation process. In our experience during computer programming, such a practice improves the calculation speed significantly than the traditional loop technique thereby allowing efficient execution of the associated codes. Furthermore, for the purpose of taking full advantage of the available computational resources, we arrange the computer codes to be amenable for execution both in the central processing unit (CPU) or in the more advanced graphical processing unit (GPU) platform.

To search out for the evasive maneuver traits exhibited by the schooling prey to avoid predation, we run several predator-prey simulations in two-dimensional space $\mathbb{R}^{2}$ under the condition of moderate prey's school size $N = 30$; the simulations are carried out for both hunting strategies. In this regard, we prescribe the simulations such that the predator only attacks the schooling prey once. Consequently, the total simulation time for all runs is fixed at $t_{\textrm{max}} = 3,000$ while the simulation time step is set to be $\Delta t = 5 \times 10^{-3}$. Furthermore, appropriate tuning upon key model parameters is undertaken as well. During the assessment, once a particular avoiding pattern is discovered, we then repeat the simulation $100$ times under the same parameter settings. Such an approach is carried out to examine the capacity of our model in returning a particular avoidance pattern consistently, considering the stochastic nature of our model. Although minor variation in the detailed structure of a particular evasive pattern is observed, the general characteristic of such a pattern is consistently obtained from all $100$ runs. Such a condition indicates that our model occupies excellent robustness in returning a particular predator-avoidance pattern.

\begin{figure}[!ht]
\centering
\includegraphics[scale=0.78]{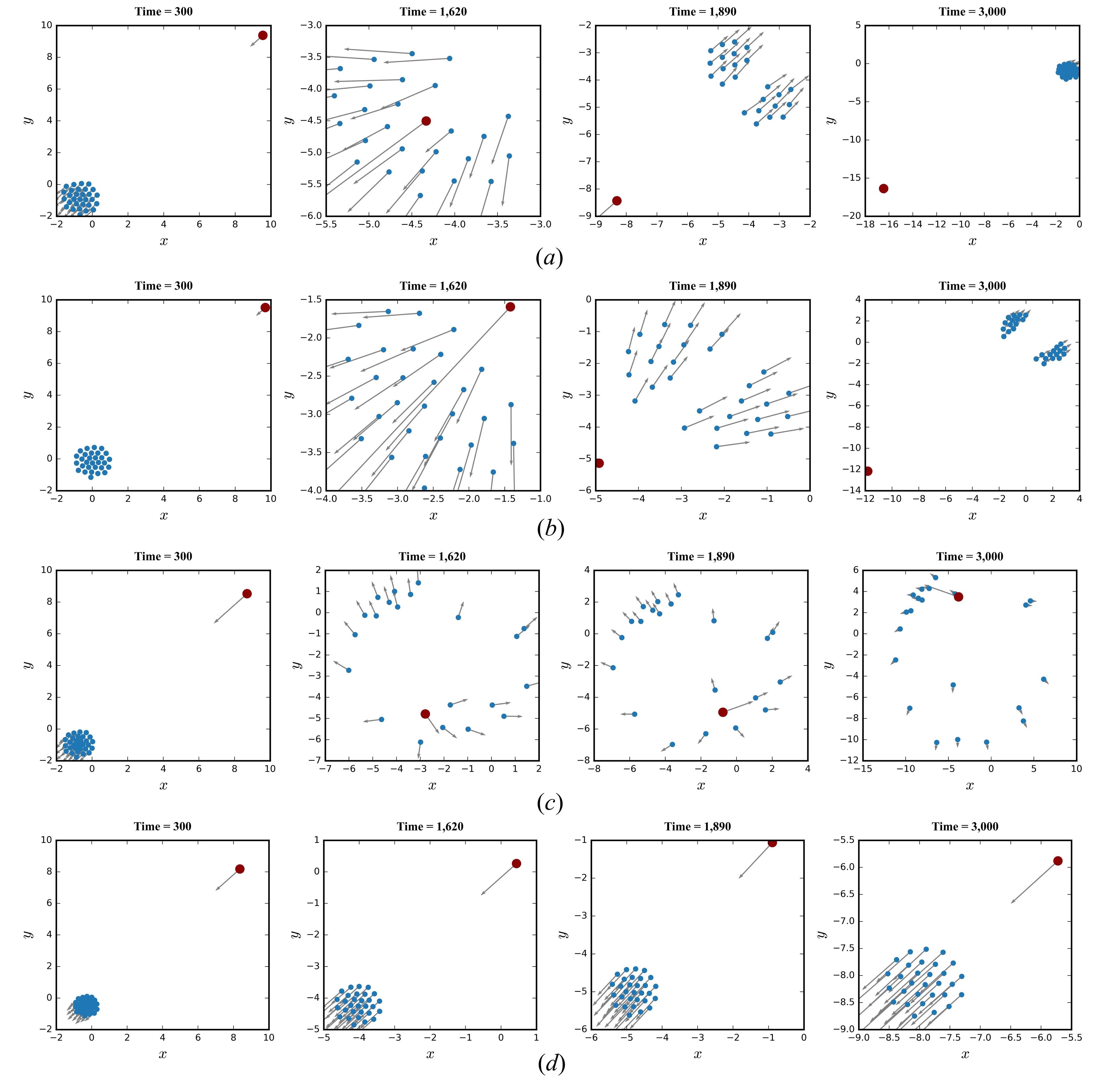}
\caption{Predator-avoidance patterns in two-dimensional space obtained from the proposed predator-prey model: ($a$) Split and Reunion, ($b$) Split and Separate into Two Groups, ($c$) Scattered, and ($d$) Maintain Formation and Distance. In the figure, the blue dots represent the schooling prey while the red dot denotes the solitary predator; the length and direction of the arrows attached to each of the dots designate the magnitude of velocity and direction of movements of each corresponding entity at the particular time, respectively.}
\label{fig3-1}
\end{figure}

Based on the results of numerical simulations, we identify four evasive patterns of the schooling prey as defensive responses to the approaching predator. We label them as: (i) Split and Reunion, (ii) Split and Separate into Two Groups, (iii) Scattered, and (iv) Maintain Formation and Distance. Figure~\ref{fig3-1} shows the corresponding predator-avoidance profiles resulted from our model. Table~\ref{table3-1} summarizes the adopted model parameter values to obtain the concomitant patterns, whereas Table~\ref{table3-2} outlines the associated predator's hunting strategy to retrieve such patterns.

\begin{table}[!ht]
\vspace*{0.2cm}
\caption{Model parameter settings for each of the discovered predator-avoidance patterns.} \label{table3-1}
\vspace*{-0.2cm}
\footnotesize \centering
\begin{tabularx}{1.02\textwidth}{ X<{\centering\arraybackslash} c c c c c c c c }
\arrayrulecolor{Black} \hline \noalign{\smallskip}
\multirow{2.5}{*}{\small Predator-Avoiding Pattern} &
\multicolumn{8}{c}{\small Model Parameter}\\[0.2cm]
& \small $\alpha$ & \small $\beta$ & \small $\delta$ & \small $p$ & \small $\theta_{1}$ & \small $\theta_{2}$ & \small $\gamma_{1}$ & \small $\gamma_{2}$\\
\noalign{\smallskip} \hline \noalign{\medskip}
Pattern I: Split and Reunion & $15$ & $0.5$ & $1$ & $4$ & $1$ & $0.5$ & $0.08$ & $0.1$\\ \noalign{\smallskip}
Pattern II: Split and Separate into Two Groups & $1$ & $0.5$ & $1$ & $4$ & $5$ & $1$ & $0.1$ & $0.1$\\ \noalign{\smallskip}
Pattern III: Scattered & $1$ & $0.5$ & $5$ & $2$ & $1$ & $2$ & $1$ & $0.1$\\ \noalign{\smallskip}
Pattern IV: Maintain Formation and Distance & $2$ & $0.5$ & $0.1$ & $2$ & $1$ & $1$ & $5$ & $10$\\
\noalign{\smallskip} \hline
\end{tabularx}
\end{table}

\begin{table}[!ht]
\vspace*{0.2cm}
\caption{\centering The adopted hunting strategy in the predator-prey simulation for \newline the attainment of predator-avoiding patterns.} \label{table3-2}
\vspace*{-0.2cm}
\footnotesize \centering
\begin{tabularx}{1\textwidth}{ m{0.55\textwidth}<{\centering} X<{\centering\arraybackslash} }
\arrayrulecolor{Black} \hline \noalign{\smallskip}
\small Predator-Avoiding Pattern & \thead{\small The Associated Hunting Strategy\\ \small Adopted in the Simulation}\\
\noalign{\smallskip} \hline \noalign{\medskip}
Pattern I: Split and Reunion & II\\ \noalign{\smallskip}
Pattern II: Split and Separate into Two Groups & I\\ \noalign{\smallskip}
Pattern III: Scattered & II\\ \noalign{\smallskip}
Pattern IV: Maintain Formation and Distance & I\\
\noalign{\smallskip} \hline
\end{tabularx}
\end{table}

Pattern I (Split and Reunion) is the first evasive maneuver trait of the schooling prey obtained from our model (see Figure~\ref{fig3-1}$a$). As can be seen in the figure, the schooling prey is initialized as a single composite unit with a circular shape in two dimensional domain, consistent with the prescribed initial condition of our numerical simulations. An observation study \citep{Larsson2012Fish} identified such a circular structure as the defensive, look-around schooling mode of the associated prey. As the predator approaches, the schooling prey reacts accordingly by exhibiting a polarized state: all individual members of the school are oriented in a parallel arrangement, facing similar direction, and move together as an egalitarian unit (the leftmost profile of Figure~\ref{fig3-1}$a$). As the predator arrived in the vicinity of the schooling prey, they readily displayed an anti-predator maneuver by expanding their formation at the right angles to the direction of predator's attack. Such a maneuver created some kind of vacuole in the middle-portion of the school (the second profile from the left of Figure~\ref{fig3-1}$a$). Following this trait, the schooling prey momentarily broke the unitary formation into two smaller groups on the left and right side of the advancing predator (the third profile from the left of Figure~\ref{fig3-1}$a$) before immediately rejoined behind the predator as a single school unit (the rightmost profile of Figure~\ref{fig3-1}$a$). It appears that this particular pattern emerged due to the strong interactions among individual prey within the school invoked in the model ($\alpha = 15$ and $p = 4$). In natural ecosystem, this particular evasive pattern closely resembles the fountain effect \citep{Partridge1982Structure, Pavlov2000PatternsSchooling, Potts1970Schooling, Shaw1978Schooling} (also called the $F$-maneuver). Such a term was introduced by \citet{Potts1970Schooling} to describe a particular anti-predation behavior exhibited by the schooling prey under medium-intensity of predation threat: the school splits into two smaller groups and rejoin behind the attacking predator \citep{Potts1970Schooling, Shaw1978Schooling}. 

Pattern II (Split and Separate into Two Groups) generally exhibits similar avoidance characteristics as Pattern I in regard to schooling prey's responses during early and mid-time predator attack (see the leftmost until the third profile from the left of Figure~\ref{fig3-1}$b$). In this particular evasive maneuver, however, after the school splits into two smaller groups, the latter formation is maintained until the end of simulation: the schooling prey does not rejoin to form a unitary structure (the rightmost profile of Figure~\ref{fig3-1}$b$). It is apparent that this particular anti-predator pattern emanated due to the weak interactions among the individual prey invoked in the simulation ($\alpha = 1$). Based on the observed similarity in behavioral traits with Pattern I, it is therefore conceivable to treat Pattern II as an adaptive form of the conventional $F$-maneuver pattern due to some kind of limitations imposed by the natural habitat of the living organisms (for example, the presence of coral reef in the vicinity that limits the evasive maneuver of the schooling prey).

In Pattern III (Scattered), the schooling prey broke the unitary formation and displayed a disordered trait. In this case, the predator's attack managed to cause total disruption of the school, rendering the individual prey to move radially outward in such an abrupt manner (see Figure~\ref{fig3-1}$c$). In regard to our model, it seems that this particular evasive maneuver appeared due to the large value of the fright factor ($\delta = 5$) imposed in the simulation. A large value of $\delta$ represents the panic condition of the prey. Such a situation is in agreement with observational studies \citep{Partridge1982Structure, Potts1970Schooling}: the schooling prey exhibits a tendency to break school formation when the school's members become panic due to the predator's attack. \citet{Potts1970Schooling} associated such a trait with a strong fright response of the schooling prey due to the high-intensity of predation threat.

The final evasive maneuver obtained from our model is Pattern IV (Maintain Formation and Distance). Here, the schooling prey demonstrates vigilant behavior: the school's members kept attention to the approaching predator and react accordingly by maintaining a safe distance from the predator. The schooling prey is readily polarized and moved in harmony with each other to avoid predation (see Figure~\ref{fig3-1}$d$). The school formation is therefore preserved until the end of the simulation period. In the framework of our model, it seems likely that this particular pattern was induced by the large value of $\gamma_{2}$ employed in the simulation ($\gamma_{2} = 10$): such a condition allows the predator to move in alignment with the schooling prey. This particular pattern shares similar traits with those of low-intensity predation threat defined by \citet{Potts1970Schooling}. 

\begin{figure}[!b]
\centering
\includegraphics[scale=0.6]{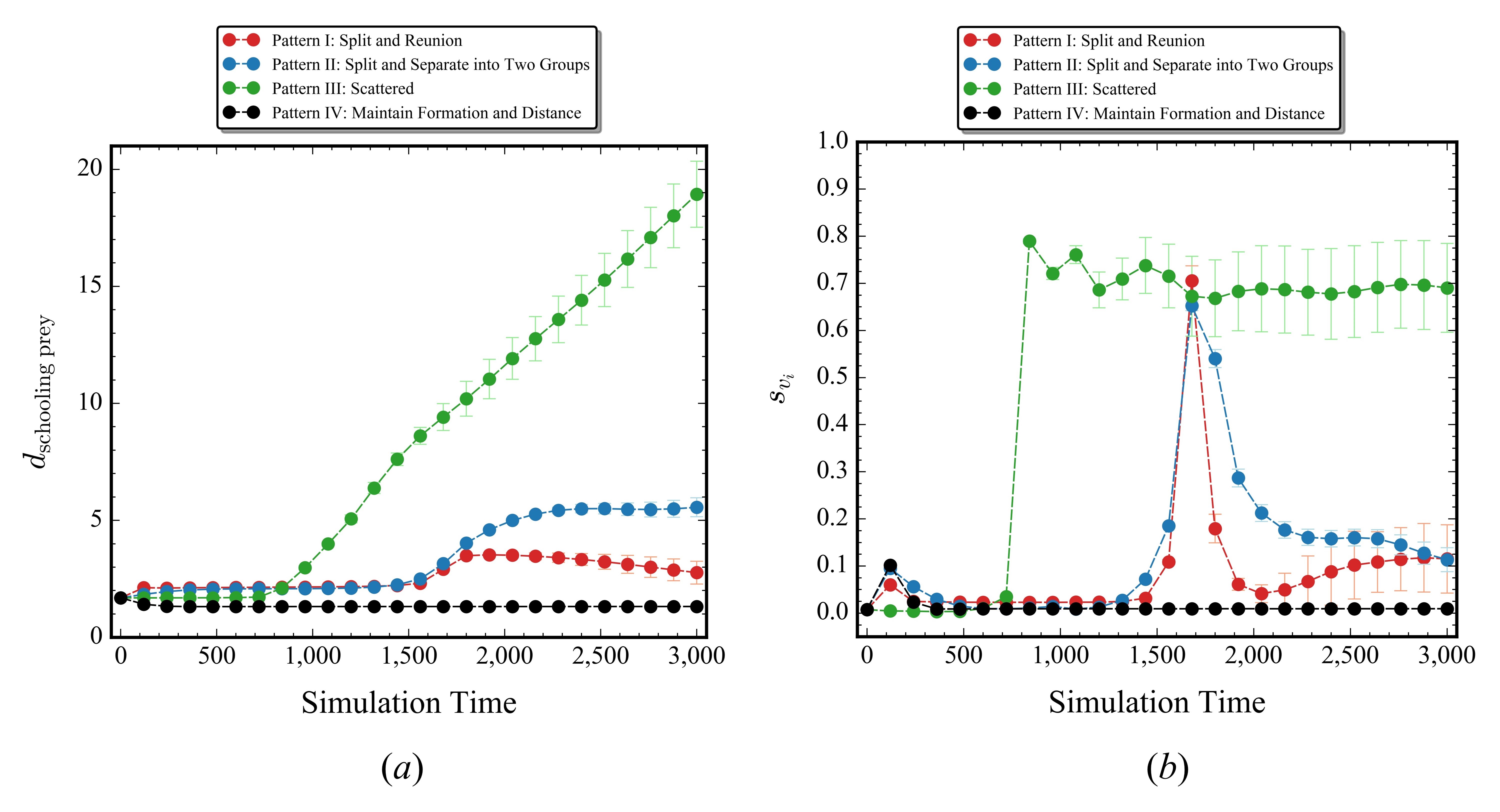}
\caption{Characteristics of the discovered predator-avoidance patterns in terms of: ($a$) Diameter of the schooling prey ($d_{\textrm{schooling prey}}$) and ($b$) Standard deviation of prey's velocity ($s_{v_{i}}$) as a function of simulation time, respectively.}
\label{fig3-2}
\end{figure}

Let us now discuss the characteristics of the discovered predator-avoidance patterns in more detail. Figure~\ref{fig3-2} displays the characteristics of the retrieved anti-predation patterns in terms of schooling prey's diameter ($d_{\textrm{schooling prey}}$) and standard deviation of prey's velocity ($s_{v_{i}}$) as a function of simulation time, respectively. Both parameters are defined accordingly as:

\begin{equation} \label{eq3-3}
d_{\textrm{schooling prey}}(t) = \sup_{1 \leqslant i \leqslant N} \|x_{i}(t) - x_{c}(t)\|, \quad 0 \leqslant t \leqslant t_{\textrm{max}},
\end{equation}

\begin{equation} \label{eq3-4}
s_{v_{i}}(t) = \sqrt{\dfrac{1}{N} \sum \limits_{i = 1}^{N} \|v_{i}(t) - v_{c}(t)\|^{2}}, \quad 0 \leqslant t \leqslant t_{\textrm{max}}.
\end{equation}

As can be seen in Figure~\ref{fig3-2}$a$, each of the discovered predator-avoidance patterns exhibits different characteristics in terms of diameter of the schooling prey. For Pattern I, as the school's formation breaks into two smaller groups, the value of $d_{\textrm{schooling prey}}$ increases. As the two smaller groups get closer, it decreases accordingly until it returns to the original value, indicating that the prey rejoins into a unitary school formation. Distinct characteristic is apparent for Pattern II where the value of $d_{\textrm{schooling prey}}$ increases gradually until the end of the simulation period, indicating that the two smaller groups of prey do not rejoin into a single formation. For Pattern III, a sharp increase of $d_{\textrm{schooling prey}}$ is observed. This condition is consistent with the fact that in the corresponding pattern, the prey breaks the school and move radially outward in an abrupt manner. The maintenance of single school formation in Pattern IV is clearly demonstrated by a constant value of $d_{\textrm{schooling prey}}$, after a small decrease at the initial stage of simulation.

The characteristics of $s_{v_{i}}$ for each corresponding evasive patterns reflect the tendency of the schooling prey (in the associated pattern) to move synchronously in the presence of predation threat. As we can see in Figure~\ref{fig3-2}$b$, Pattern I and II share relatively similar traits of $s_{v_{i}}$. Here, both patterns exhibit a sharp increase in $s_{v_{i}}$ in the middle portion of the simulation period, indicating an increase variability in the individual prey's velocity in association with the breaking of a unitary school formation into two smaller groups. After the prey splits formation and moves away from the attacking predator, they become less agitated and correspondingly reestablish the uniform characteristics at the end of the simulation period (shown by a decreasing trend of $s_{v_{i}}$ followed by reestablishment of a new \textquotedblleft stable\textquotedblright\ condition of $s_{v_{i}}$). For Pattern III, the panic condition of the prey is clearly demonstrated in the sharp increase of the corresponding $s_{v_{i}}$ at the early stage of the simulation. This is further supported by the fact that such a high $s_{v_{i}}$ value prevails until the end of the simulation. A relatively constant value of $s_{v_{i}}$ is observed for Pattern IV, indicating the maintenance of the synchronous individual movements within the school.

\begin{figure}[!hb]
\centering
\includegraphics[scale=0.6]{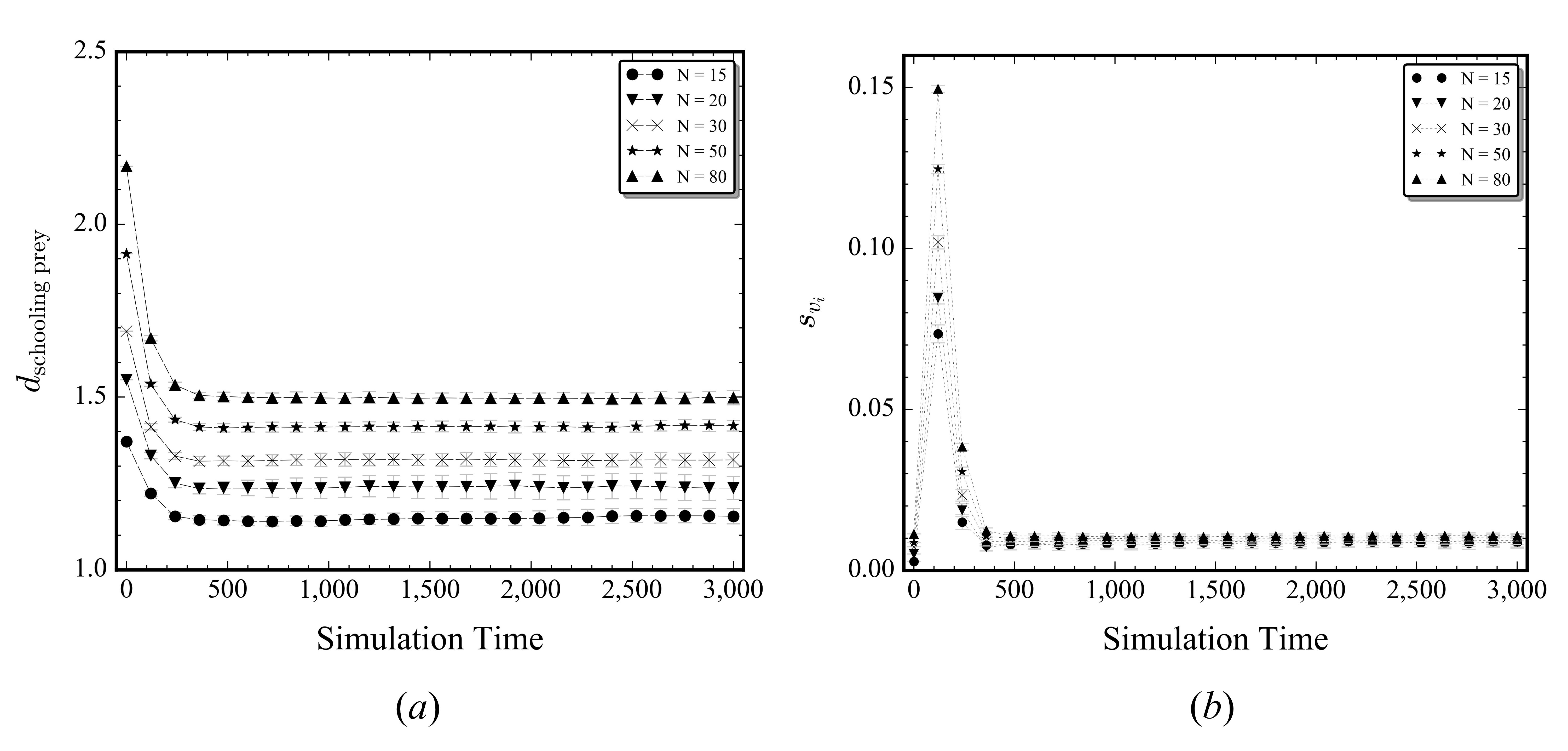}
\caption{Polarization mode displayed by the schooling prey (represented by the evasive maneuver of Pattern IV (Maintain Formation and Distance)), outlining the characteristics of: ($a$) Diameter of the schooling prey ($d_{\textrm{schooling prey}}$) and ($b$) Standard deviation of prey's velocity ($s_{v_{i}}$) as a function of simulation time, respectively.}
\label{fig3-3}
\end{figure}

We would like to recall the unique feature that occurs in Pattern IV. Therein, it is observed that the school's diameter $d_{\textrm{schooling prey}}$ decreases during the initial stage of the simulation. Such a condition reflects important characteristic of the schooling prey under imminent predation threat, namely the polarization mode of the schooling prey. \citet{Partridge1982Structure} defined such a trait as a parallel arrangement of the school members: when the school is under imminent predation threat, the individual members of the school move closer to one another and establish synchronous movements with their neighbors \citep{Partridge1982Structure}. In this regard, our model occupies excellent capability in demonstrating the occurrence of the polarization mode. To demonstrate the appearance of such a trait, we carry out numerical simulations of Pattern IV under different school size ($N = 15, 20, 30, 50, \textrm{and } 80$). Figure~\ref{fig3-3} shows the corresponding characteristics of $d_{\textrm{schooling prey}}$ and $s_{v_{i}}$, respectively. As shown in Figure~\ref{fig3-3}$a$, the first attribute of the polarization mode (namely, the individuals in the school move closer to one another) is clearly displayed by the decrease of $d_{\textrm{schooling prey}}$ during $0 \leqslant t \leqslant 500$ for all the corresponding cases. Such a condition reflects that as soon as the schooling prey realizes that they are under imminent predation threat, they immediately become more vigilant and positioned themselves closer to each other. The second attribute of the polarization mode (namely, the individuals in the school align themselves synchronously with their neighbors) is shown in Figure~\ref{fig3-3}$b$. Therein, after a sharp increase of $s_{v_{i}}$ (in association with a decrease in $d_{\textrm{schooling prey}}$) in the early stage of simulation, the value of $s_{v_{i}}$ decreases sharply and stabilizes until the end of the simulation. Such a trait applies to all the considered values of $N = 15, 20, 30, 50, \textrm{and } 80$, demonstrating consistent establishment of synchronous movements of the associated prey within the school.

\begin{figure}[!b]
\centering
\includegraphics[scale=0.55]{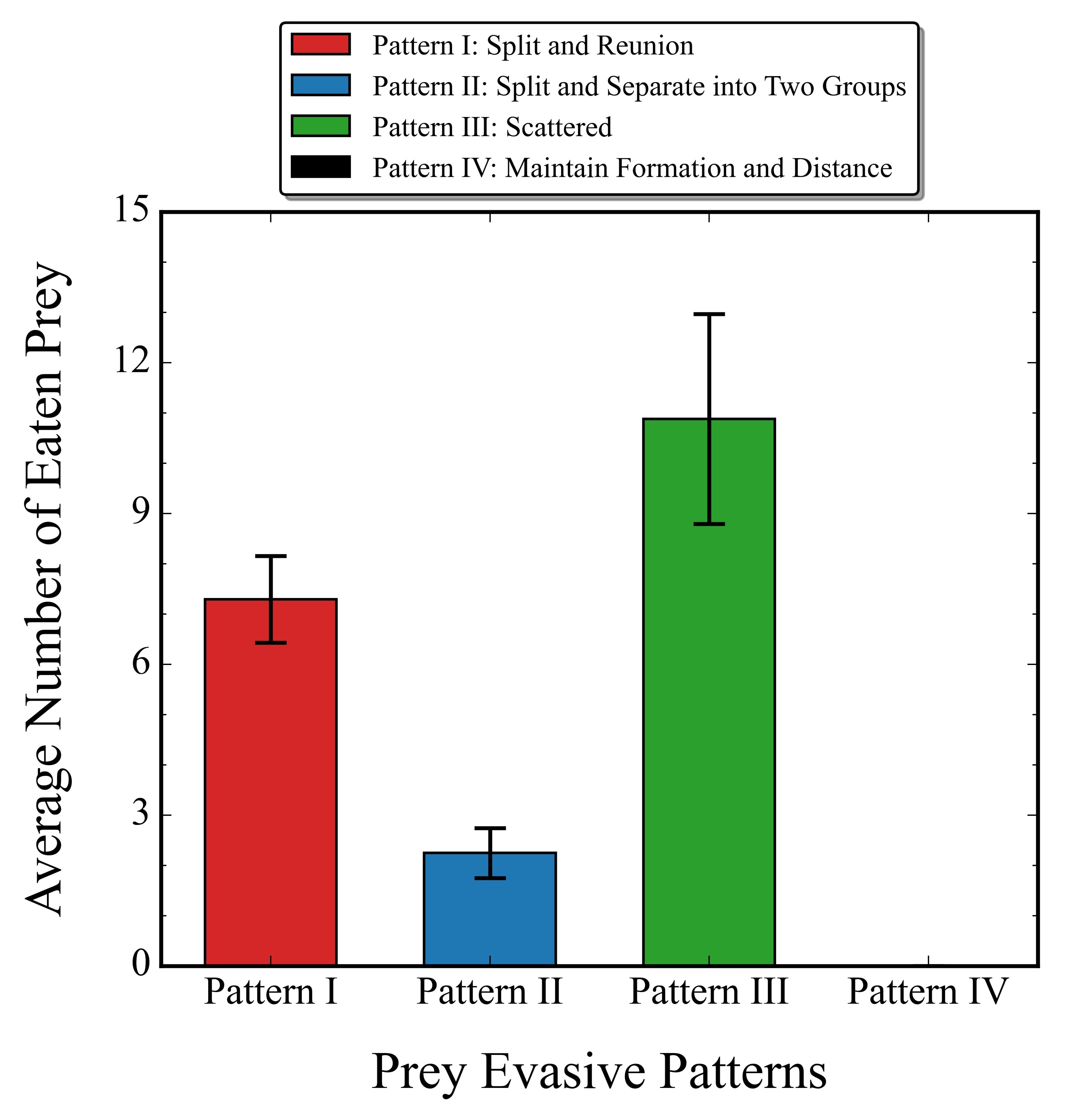}
\caption{The total number of eaten prey for each corresponding predator-avoidance patterns at the final simulation time $t_{\textrm{max} } = 3,000$.}
\label{fig3-4}
\end{figure}

Next, we would like to point out the capacity of each of the discovered predator-avoidance patterns in providing protection to the corresponding prey. To do so, we calculate the number of eaten prey at the end of the simulation period ($t_{\textrm{max}} = 3,000$) for every observed patterns. For each evasive pattern, we carry out the numerical simulations $1,000$ times. Figure~\ref{fig3-4} outlines the results of such a procedure. Therein, the largest number of eaten prey is exhibited by Pattern III (on average, up to $11$ prey is captured by the predator). From Figure~\ref{fig3-4}, it can therefore be inferred that Pattern III offers the least protection mode for the prey to avoid predation. It seems that the prey becomes more vulnerable to the predator when they are agitated and broke the unitary school formation. Such a finding is consistent with observational studies of \citet{Pavlov2000PatternsSchooling} and \citet{Shaw1978Schooling}. On the other hand, the smallest number of eaten prey is demonstrated by Pattern IV ($\textrm{the average number of captured prey} = 0$). Apparently, when the prey manages to keep the unitary school formation and maintains safe distance from the predator, the latter entity can never successfully catch an individual prey until the end of the prescribed simulation period. It is therefore conceivable to deem Pattern IV as the most effective evasive mode to avoid predation. Furthermore, our model suggests Pattern II as the second-best evasive mode of the schooling prey to mitigate predation risk (on average, only $2$ prey is captured by the predator). Pattern I, however, does not come up as the favorable mode for the schooling prey as the average number of eaten prey is observed to be four times as large as Pattern II.

As a final annotation of this section, let us point out a few relevant examples in real ecological environment related to the discovered evasive patterns from our model. By doing so, we sought to demonstrate the excellent conformity of the simulated predator-avoidance patterns with the real-life evasive patterns that occurred in the natural ecosystems. In this regard, the examples are taken mainly from predator-prey interaction in aquatic environment\textemdash including those of piscivorous birds\textemdash that include schooling behavior in the associated prey. Evasive maneuver in the form of Pattern I is well-exhibited by a school of grouper fish in avoiding the attack of the red bass predatory fish (\textit{Lutianus bohar}) \citep{Potts1970Schooling}. A similar trait is also displayed by a school of sand-eels (\textit{Ammodytes} sp.) in avoiding the incursion of mackerel \citep{Pitcher1983Predator}. On the other hand, the defensive trait of Pattern II is observed to be adopted by a school of hardyheads (\textit{Pranesus capricornensis}) in escaping the attack of avian predators, such as tern (\textit{Thalasseus bengalensis}) \citep{Hulsman2021FishyTales} and Australasian gannet (\textit{Morus serrator}) \citep{Capuska2011Dive}. Other observational studies \citep{Vabo1997Individual, Nottestad1999Herring} reported similar behavior exhibited by a school of herring (\textit{Clupea harengus}) in dodging the attack of a killer whale (\textit{Orcinus orca}). Anti-predation maneuver of Pattern III is reported to be observed in a school of Atlantic silversides (\textit{Menidia menidia}) in escaping the predation of a predatory black sea bass fish (\textit{Centropristes striata}) \citep{Pavlov2000PatternsSchooling}. Finally, the evasive maneuver of Pattern IV is found to be performed by a school of snappers (\textit{Lutianus monostigma}) in escaping the attack of a barracuda or a black-tipped shark \citep{Partridge1982Structure, Potts1970Schooling, Shaw1978Schooling}.

\section{Benefit of Prey's Schooling in Surviving Predator's Attack} \label{benefit_schooling}
\noindent In this section, we investigate the contribution of schooling prey's population size $N$ upon the chance of survival (not being eaten) of each individual prey in the school as well as their corresponding living time. 

First, let us define a mathematical expression for the probability of individual prey to be eaten by the predator:\\
\begin{equation} \label{eq4-1}
P_{\textrm{eaten}} (t) = P_{\textrm{eaten}} (t - 1) + \frac{N_{\textrm{survived}} (t - 1) - N_{\textrm{survived}} (t)}{N}.
\end{equation}
Here, $P_{\textrm{eaten}} (t)$ and $P_{\textrm{eaten}} (t - 1)$ denote the probability of individual prey being eaten at current and previous time steps, respectively while $N_{\textrm{survived}} (t)$ and $N_{\textrm{survived}} (t - 1)$ are the number of survived prey at current and previous time steps, respectively. $P_{\textrm{eaten}} (t)$ is therefore calculated recursively with initial condition set as $P_{\textrm{eaten}} (0) = 0$. The primary interest here is examining the characteristics of $P_{\textrm{eaten}}$ at the allotted final simulation time $t_{\textrm{max}}$, namely $P_{\textrm{eaten}} (t_{\textrm{max}})$.

Subsequently, the living time of each individual prey in the school ($t_{\textrm{alive}}$) is defined as the time interval during which the corresponding prey is not consumed by the predator. Such a definition implies that $t_{\textrm{alive}}$ can be determined by either: (i) tracking the time from the beginning of the simulation until the time when a prey is eaten by the predator (if the prey is eaten during the allotted simulation time) or (ii) assigning the maximum simulation time $t_{\textrm{max}}$ as the corresponding living time of the prey (if the prey is not eaten during the allotted simulation time). Thereupon, the average living time of all prey in the school ($\bar{t}_{\textrm{alive}}$) can be determined as:\\
\begin{equation} \label{eq4-2}
\bar{t}_{\textrm{alive}} =
\begin{aligned}
\begin{cases}
t_{\textrm{max}},
& \text{if $N_{\textrm{eaten}}(t_{\textrm{max}}) = 0$,}  \\
\dfrac{1}{N} \sum \limits_{i = 1}^{N} (t_{\textrm{alive}})_{i},
& \text{if $N_{\textrm{eaten}}(t_{\textrm{max}}) \neq 0$.} 
\end{cases}
\end{aligned}
\end{equation}

\noindent Here, $N_{\textrm{eaten}}(t_{\textrm{max}}) = N - N_{\textrm{survived}}(t_{\textrm{max}})$ is the number of eaten prey at final simulation time. The procedural steps for undertaking calculations of $P_{\textrm{eaten}}(t_{\textrm{max}})$ and $\bar{t}_{\textrm{alive}}$ are summarized in Scheme~\ref{algorithm_3}.

\begin{algorithm}[!ht]
\vspace*{0.2cm}
\caption{\centering Pseudocode for computing $P_{\textrm{eaten}} (t_{\textrm{max}})$ and $\bar{t}_{\textrm{alive}}$.}\label{algorithm_3}
\setstretch{1.1}
\hrulefill
\begin{algorithmic}[1]
\Function{Benefit Schooling}{$N, d, \sigma_{w}, \alpha, \beta, \delta, r, p, q, \theta_{1}, \theta_{2}, \gamma_{1}, \gamma_{2}, R_{1}, R_{2}, m, \Delta t, t_{\textrm{max}}$}
\State Load $N_{\textrm{survived}}\ (t)$ from the output of \Call{Pattern Generation}{} function;
\State Define variables $P_{\textrm{eaten}}$ and $t_{\textrm{eaten}}$; \Comment{\parbox[t]{0.48\textwidth}{$P_{\textrm{eaten}}$ \textit{is the probability of a prey to be eaten}}}\\ \Comment{\parbox[t]{0.48\textwidth}{$t_{\textrm{eaten}}$ \textit{is the time when a prey is eaten}}}
\State Define variables $t_{\textrm{alive}}$ and $\bar{t}_{\textrm{alive}}$; \Comment{\parbox[t]{0.48\textwidth}{$t_{\textrm{alive}}$ \textit{is the living time of a prey}}}\\ \Comment{\parbox[t]{0.48\textwidth}{$\bar{t}_{\textrm{alive}}$ \textit{is the average living time of all prey}}} 
\For{$t = 1$ \textbf{to} $t = t_{\textrm{max}}$}
\State Calculate $P_{\textrm{eaten}}$ using Eq.~(\ref{eq4-1});
\If{$N_{\textrm{survived}}\ (t - 1) - N_{\textrm{survived}}\ (t) > 0$}
\For{every eaten prey} 
\State $t_{\textrm{eaten}} \coloneqq t$;
\EndFor
\EndIf
\EndFor
\State Determine $N_{\textrm{eaten}}$ based on the number of members in $t_{\textrm{eaten}}$;
\For{every $i$ prey}
\If{prey $i$ was eaten}
\State $t_{\textrm{alive}} \coloneqq t_{\textrm{eaten}}$;
\ElsIf{prey $i$ was not eaten}
\State $t_{\textrm{alive}} \coloneqq t_{\textrm{max}}$;
\EndIf
\EndFor
\If{$N_{\textrm{eaten}} = 0$}
\State Calculate $\bar{t}_{\textrm{alive}}$ using the first expression in \eqref{eq4-2};
\ElsIf{$N_{\textrm{eaten}} \neq 0$}
\State Calculate $\bar{t}_{\textrm{alive}}$ using the second expression in \eqref{eq4-2};
\EndIf\\
\textbf{return} $P_{\textrm{eaten}}$, $\bar{t}_{\textrm{alive}}$
\EndFunction
\end{algorithmic}
\hrulefill
\end{algorithm}

\begin{figure}[!ht]
\centering
\includegraphics[scale=0.68]{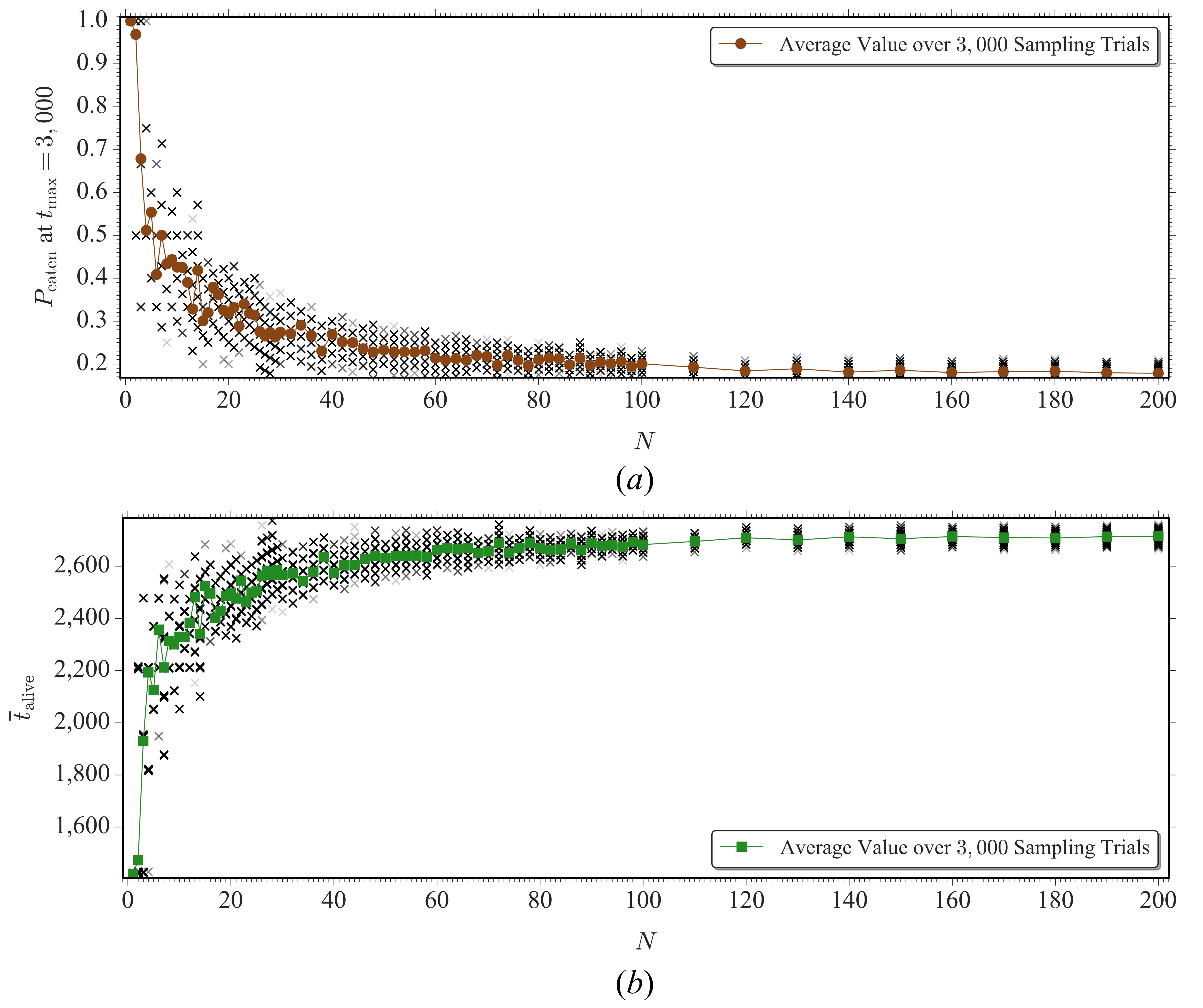}
\caption{Numerical results of center-attacking strategy: ($a$) Profiles of $P_{\textrm{eaten}} (t_{\textrm{max}})$ as a function of school size $N$, ($b$) Profiles of $\bar{t}_{\textrm{alive}}$ as a function of school size $N$.}
\label{fig4-1}
\end{figure}

\begin{figure}[!ht]
\centering
\includegraphics[scale=0.68]{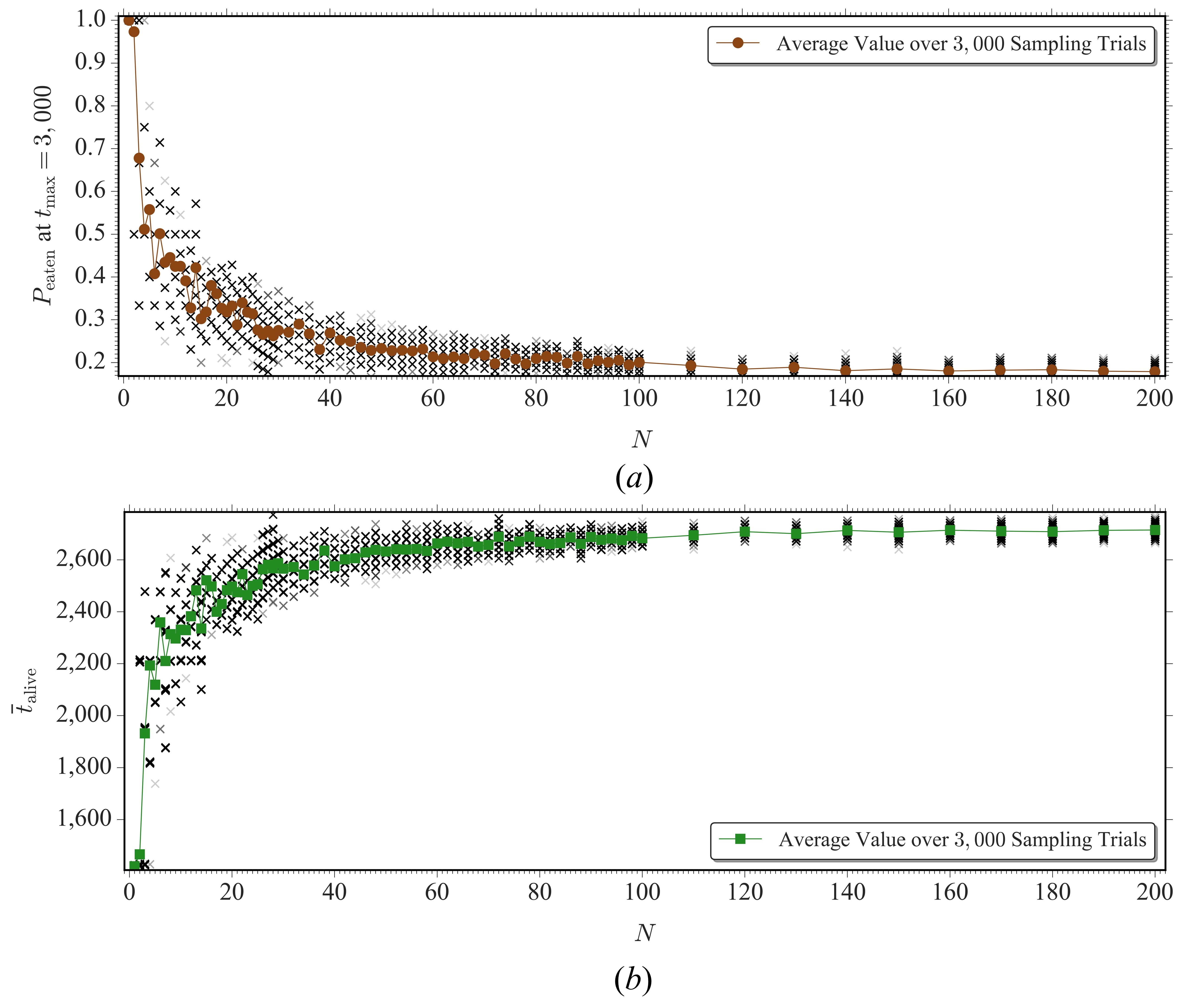}
\caption{Numerical results of nearest-attacking strategy: ($a$) Profiles of $P_{\textrm{eaten}} (t_{\textrm{max}})$ as a function of school size $N$, ($b$) Profiles of $\bar{t}_{\textrm{alive}}$ as a function of school size $N$.}
\label{fig4-2}
\end{figure}

To understand the influence of prey's population size $N$ upon the characteristics of $P_{\textrm{eaten}}(t_{\textrm{max}})$ and $\bar{t}_{\textrm{alive}}$, numerical simulations of Schemes~\ref{algorithm_2} and \ref{algorithm_3} are executed for different values of $N$. In this regard, key simulation parameters are kept constant for all simulations at $\alpha = 15$, $\beta = 1$, $\delta = 1$, $p = 4$, $\theta_{1} = 0.1$, $\theta_{2} = 0.5$, $\gamma_{1} = 0.1$, and $\gamma_{2} = 0.1$. To ensure that the predator only attacks the schooling prey once, we maintain the total simulation time at $t_{\textrm{max}} = 3,000$. The same value of simulation time step $\Delta t = 5 \times 10^{-3}$ is prescribed, as well. As part of an effort to see the benefit of schooling in escaping predation, we run the simulations gradually starting from solitary prey ($N = 1$) up to sufficiently large value of prey's population ($N = 200$). Furthermore, to average out the contribution of the stochastic effect in our model, predator-prey simulation under one condition of $N$ is executed repeatedly up to a sufficiently large number of repetitions. Here, the repetitions (sampling trials) are carried out $3,000$ times so as to obtain adequate data to capture inherent variability induced by the stochastic nature of our model. Clearly, such an operation involves a large number of simulations. We therefore implement a parallel computational process. Such a configuration allows the simulations to be executed up to ten times faster than the conventional serial computing.

We carry out the above procedures for the two proposed hunting strategies. Figure~\ref{fig4-1} shows the behavior of $P_{\textrm{eaten}} (t_{\textrm{max}})$ and $\bar{t}_{\textrm{alive}}$ for the center-attacking strategy, while Figure~\ref{fig4-2} demonstrates similar behavior for the nearest-attacking tactic.

Let us first discuss the characteristics of $P_{\textrm{eaten}} (t_{\textrm{max}})$ as a function of school size $N$. As can be seen in Figures~\ref{fig4-1}(a) and \ref{fig4-2}(a), in general, the values of $P_{\textrm{eaten}} (t_{\textrm{max}})$ for both hunting strategies decrease as $N$ increases. These values can be divided into three general trends. The first trend occurs at $1 \leqslant N \leqslant 20$. During such a period, a sharp decrease of $P_{\textrm{eaten}} (t_{\textrm{max}})$ is observed for both attacking tactics, accompanied by relatively wide variability among the sampling data. A clear benefit of schooling in providing better safety from predation is displayed during this particular period: the probability of individual prey to be eaten by the predator reduced threefold when the prey constitutes a small group ($N = 20$) than when it is wandering alone ($N = 1$). Solitary prey is always eaten by the predator ($P_{\textrm{eaten}} (t_{\textrm{max}})$ = 1). As the function $P_{\textrm{eaten}} (t_{\textrm{max}})$ proceeds upon the second segment ($20 \leqslant N \leqslant 120$), it declines further. However, the decline rate is not as sharp as in the first segment. At $N = 120$, the value of $P_{\textrm{eaten}} (t_{\textrm{max}})$ is approximately 1.5 times lower than its value at $N = 20$. In this second trend, the amount of variability due to the stochastic effect is found to decrease as well. Finally, at the third segment ($120 \leqslant N \leqslant 200$), the contribution of $N$ becoming stable and further increase in $N$ does not reduce the probability of being eaten any longer. It is therefore reasonable to infer that at this stage, the system has reached its saturation point in terms of schooling advantage in avoiding predation. Although minor variability is observed, the general trend is almost stable at $P_{\textrm{eaten}} (t_{\textrm{max}}) \approx 18\%$. Thus, our model suggests that schooling among prey is capable of reducing the probability of individual prey to be eaten from 100\% (solitary condition) up to around 18\% (by constituting a school of $N = 200$).

The characteristic of $\bar{t}_{\textrm{alive}}$ against school size $N$ demonstrates reciprocal behavior with the trend of $P_{\textrm{eaten}} (t_{\textrm{max}})$. The general trait is, however, equivalent with the latter quantity. Lower predation risk associated with increasing $N$ results in longer average living time $\bar{t}_{\textrm{alive}}$ of the schooling prey (see Figures~\ref{fig4-1}(b) and \ref{fig4-2}(b)). 

At this point, we can extract one important characteristic from the outcomes of numerical simulations: as the size of the schooling prey increases, the probability of individual prey to be eaten by the predator decreases. We find that such a trait is in excellent agreement with empirical evidence. \citet{Radakov1958Adaptive} reported experimental investigation involving school of coalfish prey (\textit{Pallachius virens}) under predation threat of cod (\textit{Gadus morhua}). Therein, it was observed that when the coalfish existed as a solitary individual, the cod took only 26 seconds to capture it. The corresponding hunting time of the cod was found to increase five-fold when the prey coalfish resides in a school of 25 to 35 individuals, demonstrating the advantage of prey schooling. In a separate study, \citet{Neill1974Experiments} designed a controlled experiment involving several predatory aquatic animals (squid, cuttlefish, pike, and perch) alongside their associated prey (mullet, bleak, dace, guppy, and minnow). They found that for all predator-prey systems, increasing the number of prey from 1 to 6 and then to 20 reduced the probability of success of the predator in capturing the prey. The functionality of prey group size in reducing the risk of predation is also reported in \citep{Partridge1982Structure, Shaw1978Schooling, Parrish1989Reexamine}.

As a final annotation, we would like to relate our findings with the established theoretical framework of predator-prey interaction. In this regard, mitigation of predation risk due to increasing group size in the schooling prey species has been acknowledged as the dilution effect \citep{Ioannou2021Grouping, Pitcher1986Functions}. \citet{Pitcher1986Functions} argued that such an effect should be considered alongside the number of predator's attack upon the associated schooling prey, resulting in the so-called attack-abatement effect. In our simulation study, we realize that the advantage of schooling ceases as the prey's population size $N$ reaches the third segment ($120 \leqslant N \leqslant 200$): as $N$ increases above $120$, there is no additional reduction in the probability of individual prey being eaten by the predator. Such a condition may be attributed to the fact that we limit the number of predator's attack in our simulation: only single attack of the predator is allowed. Therefore, in the light of theoretical framework, our model occupies excellent capability in demonstrating the dilution effect. Consideration of the attack-abatement effect is, however, not within the scope of this study. Such a point of view is left for future work.

\section{Conclusions} \label{conclusions}
\noindent In this study, we propose mathematical model to describe predator-prey system with schooling behavior imposed upon the prey. The model is of stochastic differential equations (SDEs), where interactions between the predator and the schooling prey are manifested through attractive or repulsive forces analogous to those occurred in molecular physics. In addition, two hunting tactics of the predator are proposed and integrated into the general model, namely the center-attacking and the nearest-attacking strategy. We carry out numerical simulations to test the capability of our model in returning predator-avoidance patterns typically exhibited by the schooling prey as defensive responses to predation threat. Furthermore, we also address the capacity of our model in showing the benefit of constituting large school of prey in mitigating the risk of predation of an individual prey. The key findings are as follows:
\begin{enumerate}[label={(\arabic*)}]
\itemsep-0.2cm
\item On the basis of numerical simulations upon our model, four predator-avoidance patterns are discovered. We call them: (i) Split and Reunion, (ii) Split and Separate into Two Groups, (iii) Scattered, and (iv) Maintain Formation and Distance. These simulated patterns are in agreement with real-life evasive patterns exhibited by diverse species of prey fish in their natural habitat. 
\item In the framework of the discovered predator-avoidance patterns, the last evasive pattern (Maintain Formation and Distance) is observed to yield the best protection mode for the schooling prey in escaping predation threat. Meanwhile, the third evasive pattern (Scattered)\textemdash which manifests panic condition of the schooling prey\textemdash is found to provide the least protection mode. Such findings are consistent with observational viewpoint.
\item Our model successfully demonstrates the benefit of constituting large prey school in mitigating predation risk. In general, it is observed that as the prey's population in the school increases, the probability of an individual prey to be eaten by the predator decreases. This result is consistent with observations in the real-life aquatic ecosystem.
\item Our model is consistent with the established theoretical framework of predator-prey system: the model possesses excellent capability in demonstrating the occurrence of the polarization mode of the schooling prey and the dilution effect.
\end{enumerate}


\section*{Acknowledgment}
\noindent
The work of the first (A.D.H.) and last (T.V.T.) authors were supported by JSPS KAKENHI Grant Number 19K14555.


\bibliographystyle{elsarticle-num-names-nodoiprefix}
\bibliography{mybib}

\end{document}